\documentclass[12pt]{article}
\usepackage{amsmath,amstext,amsthm,amssymb,amsfonts,latexsym,amscd,bezier}
\usepackage[latin1]{inputenc}
\usepackage{color}
\usepackage{hyperref}
\usepackage{pdfsync}
\usepackage[all]{xy}
\usepackage{graphicx}
\newtheorem{theorem}{Theorem}

\newtheorem{remark}{Remark}

\newcommand{\bb}{\mathbb}
\newcommand{\R}{\bb R}
\newcommand{\C}{\bb C}
\newcommand{\N}{\bb N}
\newcommand{\s}{\Bbb{S}}
\renewcommand{\o}{\overline}

\newcommand{\n}{\noindent}

\renewcommand{\phi}{\varphi}
\renewcommand{\epsilon}{\varepsilon}

\newcommand{\Img}{\operatorname{Im}}
\newcommand{\Ker}{\operatorname{Ker}}
\newcommand{\coker}{\operatorname{coker}}

\begin{document}

\title{The $C^r$ dependence problem of eigenvalues of the Laplace operator on domains in the plane
}

\author{\textbf{Julian Haddad \footnote{\textit{E-mail addresses}:
julianhaddad@ufmg.br (J. Haddad)}}\\ {\small\it Departamento de
Matem\'{a}tica, Universidade Federal de Minas Gerais,}\\ {\small\it
Caixa Postal 702, 30123-970, Belo Horizonte, MG, Brazil}\\
\textbf{Marcos Montenegro \footnote{\textit{E-mail addresses}:
montene@mat.ufmg.br (M. Montenegro)}}\\ {\small\it Departamento de
Matem\'{a}tica, Universidade Federal de Minas Gerais,}\\ {\small\it
Caixa Postal 702, 30123-970, Belo Horizonte, MG, Brazil}} \maketitle

\markboth{abstract}{abstract}
\addcontentsline{toc}{chapter}{abstract}

\hrule \vspace{0,2cm}

\n {\bf Abstract}

The $C^r$ dependence problem of multiple Dirichlet eigenvalues on domains is discussed for elliptic operators by regarding smooth one-parameter families of $C^1$ perturbations of domains in $\R^n$. As applications of our main theorem (Theorem 1), we provide a fairly complete description for all eigenvalues of the Laplace operator on disks and squares in $\R^2$ and also for its second eigenvalue on balls in $\R^n$ for any $n \geq 3$. The central tool used in our proof is a degenerate implicit function theorem on Banach spaces of independent interest.

\vspace{0.5cm}
\hrule\vspace{0.2cm}

\section{Introduction and statement}

A number of eigenvalue problems associated to the Laplace operator have been widely investigated since the famous Rayleigh's book \cite{R}. We refer for instance to the Henrot's \cite{He} book for details on background material.

The present paper deals with the dependence problem of Dirichlet eigenvalues for elliptic operators with respect to perturbations of the domain of embedding type. This is a classical problem in the literature which has been addressed at different times, mainly regarding algebraically simple eigenvalues, namely eigenvalues with algebraic multiplicity equal to $1$. For this case we refer among other works to \cite{R}, \cite{Had}, \cite{CH}, \cite{BV}, \cite{CL}, \cite{D}, \cite{C}, \cite{P}, \cite{H} and \cite{HM}.

In the case that the eigenvalues are algebraically multiple, much less results on dependence are known. Indeed, one does not hope in general that multiple eigenvalues depend differentiably on smooth perturbations of the domain. Examples can easily be constructed even in finite dimension (see for example page 37 of \cite{He}). Nonetheless, weaker differentiability of eigenvalues with respect to differentiable deformations of $C^2$ domains has been investigated in the multiple situation by Cox, namely in the sub differentiability (or Lipschitz) sense (see Theorem 1 of \cite{C}), and by Rousselet \cite{Ro} and Munnier \cite{Mu}, in the lateral directional derivative sense (see Theorem 2.5.8 of \cite{He}). On the other hand, the $C^r$ dependence of multiple eigenvalues upon smooth one-parameter families of perturbations of the domain seems to be unknown even for $r = 1$ and smooth perturbations. Important contributions can also be found in the classical Kato's book \cite{K} where are provided expansions of eigenvalues considering one-parameter perturbations in the analytical context. For an excellent overview on dependence of eigenvalues with respect to the domain, among other interesting problems, we also refer to the Henrot's \cite{He} book.

According to the works of Uhlenbeck \cite{U1, U2} and Pereira \cite{Per}, it is interesting to note that algebraic simplicity of eigenvalues of the Laplace operator is a generic property. Precisely, given a bounded open subset $\Omega$ of $\R^n$ and a number $\varepsilon > 0$, there exists a diffeomorphism $\varphi: \Omega \to \varphi(\Omega) \subset \R^n$ such that $\|\varphi - Id\|_{C^3} \leq \varepsilon$ and $\varphi(\Omega)$ has the following property: all eigenvalues of the Laplace operator on $\phi(\Omega)$ are algebraically simple. Nevertheless, eigenvalues of symmetric domains ({\it e.g.} balls and cubes in $\R^n$), except the first one, have in general algebraic multiplicity greater than $1$ (see \cite{He}).

More particularly and surprisingly, the following question is open for any $r \geq 1$:

\begin{center}
Does the Dirichlet spectrum of the Laplace operator on balls in $\R^n$ vary $C^r$ smoothly upon one-parameter perturbations of $C^1$ class of the domain for any dimension $n \geq 2$?
\end{center}

By using an appropriate degenerate implicit function theorem on Banach spaces, we establish a result on $C^r$ dependence of Dirichlet eigenvalues of elliptic operators in the multiple case (Theorem 1). As a byproduct, we give a fairly complete answer, in an almost everywhere sense, to the problem stated above in dimension $n = 2$ on disks and squares and also for the second eigenvalue on balls in $\R^n$ for any dimension $n \geq 3$. As we shall see, perturbations of the domain (which can be a certain non-smooth domain) mean images by diffeomorphic $C^1$ maps and deformations of these perturbations are considered $C^{r+1}$ with respect to the parameter.

Before we go further and state our main theorem, a suitable framework of the problem of interest should first be introduced.

The eigenvalue problem for the Laplace operator under Dirichlet boundary condition is given by

\begin{equation} \label{EP}
\left\{
\begin{array}{rcll}
- \Delta u &=& \lambda_{\Omega} u & {\rm in} \ \ \Omega,\\
u &=& 0  & {\rm on} \ \ \partial\Omega,
\end{array}
\right.
\end{equation}

\n where $\Omega$ denotes a bounded open subset of $\R^n$, not necessarily smooth, $n \geq 2$ and $\lambda_{\Omega} \in \R$ is a eigenvalue of (\ref{EP}).

Let $\Omega_0$ be a bounded open subset of $\R^n$. Roughly, let $\Omega_t$ be a smooth one-parameter perturbation of $\Omega_0$ around $t = 0$. We are interested in studying the $C^r$ regularity of the map $t \mapsto \lambda_{\Omega_t}$ around $t = 0$ in the case that $\lambda_{\Omega_0}$ is an algebraically multiple eigenvalue. The techniques often used in the context of algebraically simple eigenvalues base on bifurcation theory, implicit function theorem, perturbation theory of operators and transversality theorem. Such tools do not apply directly to the multiple case, so instead we approach the problem through a degenerate implicit function theorem on Banach spaces to be proved in the Section 2 (Theorem \ref{DIft}).
The idea is simple and consist in constructing a suitable map with invertible derivative by using the first and second derivatives of the map for which we seek an implicit function.
The method may vaguely resemble a Lyapunov-Schmidt reduction process (to finite dimension) but no Fredholm condition is required and also there is no additional nonlinear equation in our assumptions. These characteristics give Theorem \ref{DIft} an independent interest and possibly other interesting applications.

We now make precise the meaning of the expression ``smooth one-parameter perturbation of $\Omega_0$" based on the Uhlenbeck's idea \cite{U2} of parameterizing domains as images of a fixed domain via diffeomorphisms. The advantage is that the collection of such diffeomorphisms is a subset of a certain normed vector space and, therefore, has differentiable structure.

Let $\Omega_0$ be a bounded open subset of $\R^n$. Denote by $X$ the Banach space of $C^1$ maps $\phi : \Omega_0 \rightarrow \R^n$, which
extend continuously up to the boundary of $\Omega_0$ as well as its
derivatives of first order, endowed with the usual norm

\[
||\phi||_X := \max_{x \in \overline{\Omega}_0} \{ \|\phi(x)\|, \|D\phi(x)\| \}\, .
\]

\n Let $E^1(\overline{\Omega}_0)$ be the collection of maps $\phi \in
X$ such that $\varphi$ is a diffeomorphism on its image. Note that the inclusion
function $1_{\Omega_0}$ belongs to $E^1(\overline{\Omega}_0)$.

For a bounded open subset $\Omega$ of $\R^n$, we denote by $\Lambda(-\Delta, \Omega) \subset \R \times H^1_0(\Omega)$ the set of couples $(\lambda_{\Omega}, u)$ satisfying
\eqref{EP} where $u$ is nonzero.

Our main theorem states that

\begin{theorem} \label{teorema1}
Let $\Omega_0$ be a bounded open subset of $\R^n$ satisfying the divergence theorem for an unit normal vector field $\nu$ defined on $\partial \Omega_0$ in almost everywhere and oriented outwards $\Omega_0$. Let $\varphi: \R \to E^1(\overline{\Omega}_0)$ be a $C^{r+1}$ curve, $r \geq 1$, such that $\varphi(0) = 1_{\Omega_0}$. Denote $\overline \varphi = \varphi'(0)$. Let $\lambda_0$ be a multiple eigenvalue of the Laplace operator on $\Omega_0$. Assume that its eigenspace is generated by a basis $\{u_0, u_1, \ldots, u_k\} \subset H^1_0(\Omega_0)$ with $k \geq 1$.
Consider the $(k+1) \times (k+1)$ real matrices $A$ and $B$ with elements, respectively,

\[
A_{ij} = (u_i, u_j)_{\overline\varphi}: = \int_{\partial \Omega_0} \frac{\partial
u_i}{\partial \nu} \frac{\partial u_j}{\partial \nu}\; \o \phi \cdot \nu\; dS
\]

\n and

\[
B_{ij} = \langle u_i, u_j \rangle_{L^2} := \int_{\Omega_0} u_i u_j\; dx\, .
\]
Assume that the polynomial $\chi(s) = \det(A - sB)$ has a simple real zero $\mu$. Assume also that the eigenfunctions $u_i$ are of $C^1$ class up to the boundary. Then, there exist $C^{r}$ functions $t \in (-\varepsilon, \varepsilon) \mapsto \lambda(t)$ and $t \in (-\varepsilon, \varepsilon) \mapsto u(t)$ such that $\lambda(0) = \lambda_0$, $u(0) \in \langle u_0, u_1, \ldots, u_k \rangle$ and $(\lambda(t), u(t)) \in \Lambda(-\Delta, \varphi_t(\Omega_0))$ for $t \in (-\varepsilon, \varepsilon)$. Moreover, we have $\lambda'(0) = -\mu$.
\end{theorem}

In the above theorem, the assumption that eigenfunctions corresponding to $\lambda_0$ belong to $C^1(\o\Omega_0)$ can be relaxed by assuming a $L^2$ integrability of trace of their gradient. This is a case for a wide class of non-smooth domains.

The matrix $A$ has already appeared in Theorem 3.2 of \cite{Ro} (see also Theorem 2.5.8 of \cite{He}) when the basis $\{u_0, u_1, \ldots, u_k\}$ is $L^2$-orthonormal, where its eigenvalues are proved to be the possible values of the lateral directional derivative of $\lambda(t)$ at $t = 0$. The statement of Theorem \ref{teorema1} requires the zero $\mu$ of $\chi(s)$ to be simple so that no other branch of eigenvalues $\lambda(t)$ has the same derivative at $0$. This key assumption allow us to prove the $C^r$ smoothness of $\lambda(t)$ around $t = 0$ with the aid of the perturbative tool of the next section.

The paper is organized into three sections. Section 2 is devoted to the proof of our main tool, the degenerate implicit function theorem on Banach spaces. In Section 3 we present an extension of Theorem 1 to more general elliptic operators and provide its proof. Finally, in Section 4 we solve the proposed problem on disks and squares and also for the second eigenvalue on balls in $\R^n$ for almost every perturbation $\varphi$ of the domain.

\section{An abstract perturbative theorem}

This section is devoted to the proof of the following degenerate implicit function theorem.
For the concepts of infinite dimensional submanifold and other definitions from infinite dimensional differential geometry, we refer to \cite{MR} or \cite{La}.

\begin{theorem} \label{DIft} Let $X$ and $Y$ be Banach spaces and $F:\R \times X \to Y$ be a $C^{r+1}$ map, $r \geq 1$, and $M \subset X$ be a $C^{r+1}$ submanifold with complemented tangent space. Let $p$ be a fixed point on $M$ and set $q = F(0, p)$ and $L = F_x(0,p) := D_XF(0,p)$. Assume the following conditions are satisfied:
\begin{itemize}
\item[(a)] $F(0, x) = q$ for all $x \in M$;

\item[(b)] $\Ker(L) = T_p M$;

\item[(c)] $F_t(0,p) := D_{\R}F(0, p) = L(v)$ for some $v \in X$;

\item[(d)] The composition
	\begin{displaymath}
		\xymatrix
		{
			T_p M \ar[r]^-G & Y \ar[r]^-{\pi} &\coker(L)
		}
	\end{displaymath}
\n is a vector space isomorphism, where $G = D(F_t - F_v) (0,p)$ and $\pi$ is the canonical projection.
\end{itemize}

\n Then, there exists a $C^{r}$ curve $x: (-\varepsilon, \varepsilon) \to X$ defined in a neighborhood of $0$ such that $x(0) = p$ and $F(t,x(t)) = q$ for all $t \in (-\varepsilon, \varepsilon)$. Moreover, we have $x'(0) + v \in T_p M$.
\end{theorem}

\begin{remark}
The co-kernel defined as the quotient $\coker(L) = \frac Y {\Img(L)}$ needs not in general to be a Banach space but just a vector space. In our applications, $\Img(L)$ is finite-codimentional and thus, complemented.
\end{remark}

\n {\bf Proof of Theorem 2.} By direct computation, it is easy to see that all properties are preserved by $C^{r+1}$ changes of coordinates. In other words, if $\Psi: X \to \tilde X$ is a local $C^{r+1}$ diffeomorphism, then the function $\tilde F(t, \tilde x) = f(t,\Psi^{-1}(\tilde x))$ also satisfies the conditions (a)-(d). In fact, if $\tilde M = \Psi(M)$, $\tilde p = \Psi(p)$, $\tilde q = q$, $A = D\Psi(p), \tilde L = D \tilde F(0,\tilde p), \tilde v = A(v)$ and $\tilde G = D(\tilde F_t - \tilde F_{\tilde v})(0,\tilde p)$, then $\tilde L \circ A = L$ and $\tilde G \circ  A = G + L\circ D(\Psi^{-1}_{\tilde v})(\tilde p) \circ A$, so that the claim follows.

By taking a suitable coordinate neighborhood at $p \in U \subset X$ adapted to $M$, we may assume that there is a decomposition $X = X_1 \oplus X_2$ such that $M \cap U  = (\{0\} \times X_2) \cap U$. Moreover, we assume that $q=0$. By the assumptions (a) and (c), $F: \R \times X_1 \times X_2 \to Y$ satisfies
	
\[
F(0,0,x_2) = 0
\]

\n for all $x_2 \in X_2$ around $0$ and

\[
F_t(0,0,0) = L(v)\, .
\]

\n In addition, by (b), the partial derivative of $F$ in the $X_1$-direction at $(0, 0, 0)$, denoted by $L^{X_1}: X_1 \to Y$, is injective since $Ker(L) = \{0\} \times X_2$.

Let $w \in X_1$ be the first coordinate of $v$, so $L^{X_1}(w) = L(v)$. Now consider the functions
	
\[
H(t,x_1,x_2) = F(t, t(x_1 - w), x_2)
\]
and
	
\[
J(t, x_1, x_2) = \frac 1 t H(t,x_1,x_2) = \int_0^1 H_t(s \,t, x_1, x_2)\; ds
\]
which is an integral of $C^{r}$ functions.
By Example 2.4.16 (Differentiating Under the Integral) in \cite{MR}, the function $J$ is $C^{r}$.
%A standard computation gives that $J$ is also of $C^{r-1}$ class.
%XXXXXXXXXXXXXXXXXXXXXXX
%
%no trivial!!
%
%aca hay que usar que $X$ es compactamente generado
%si $h$ es una derivada parcial y $t_n, x_m \to t,x$ entonces
%$h$ es uniformemente continua en $[0,1] \times \{x, x_1, x_2, \ldots \}$
%porque es compacto. Luego $h(.,x_n)$ converge unifirmemente a $h(., x)$
%
%XXXXXXXXXXXXXXXXXXXXXXX

On the other hand,

\[
J(0,x_1,x_2) = H_t(0,x_1,x_2) = DF(0,0,x_2)[1, x_1 - w, 0]
\]
and
	
\[
DJ(0,0,0)[0, \bar x_1, \bar x_2] = DF(0,0,0)[0,\bar x_1,0] + D^2 F(0,0,0)[(1,-w,0),(0,0,\bar x_2)]
\]
	
\[
=: L^{X_1}[\bar x_1] + G[\bar x_2]\, .
\]

Note that the bounded operator $L^{X_1}[\bar x_1] + G[\bar x_2]$ is a vector space isomorphism by the assumption (d), thus a Banach space isomorphism. In addition, $J(0, 0, 0) = DF(0, 0, 0)[1, -w, 0] = F_t(0, 0, 0) - L^{X_1}[w] = 0$. So, the implicit function theorem applied to $J$ provides $\varepsilon > 0$ and $C^{r}$ functions $\tilde x_1(t)$ and $\tilde x_2(t)$ such that

\[
0 = J(t, \tilde x_1(t),  \tilde x_2(t)) = \frac 1 t F(t, t (\tilde x_1(t) - w), \tilde x_2(t))
\]
for all $t \in (-\varepsilon, \varepsilon)$ with $t \neq 0$. Therefore, the result follows by taking $x(t) = (x_1(t), x_2(t))$, where $x_1(t) = t (\tilde x_1(t) - w)$ and $x_2(t) = \tilde x_2(t)$. Indeed, the above equality yields $F(t, x(t)) = 0$ for all $t \in (-\varepsilon, \varepsilon)$. Moreover, since $x_1'(0) = -w$, we have $x'(0) + v \in \{0\} \times X_2 = T_p M$.\hfill$\Box$

\section{An extension of Theorem 1 and its proof}

In this section we state a version of Theorem 1 extended to non-self-adjoint linear elliptic operators of second order and present its proof.

For that purpose, consider a fixed second order elliptic operator ${\cal L}_0$ under the form

\[
{\cal L}_0 := \sum_{i, j = 1}^n \partial_i (a^0_{ij}(x) \partial_j) - \sum_{i = 1}^n b^0_i(x) \partial_i - c^0(x)
\]

\n with coefficients $a^0_{ij}$, $b^0_i$ and $c^0$ in $C^{r+1}(\R^n)$, $r \geq 0$. Note that any second
order linear differential operator with coefficients in $C^1(\R^n)$
can always be placed into the above format.

Let $\Omega_0$ be a bounded open subset of $\R^n$ and $X$ and $E^1(\overline{\Omega}_0)$ be as defined in the introduction. Given a bounded open subset $\Omega$ of $\R^n$, we denote by
$\Lambda({\cal L}_{0},\Omega) \subset \R \times H^1_0(\Omega)$ the set of pairs $(\lambda_{\Omega}, u)$, where $u$ is nonzero, satisfying

\begin{equation} \label{EP1}
\left\{
\begin{array}{rcll}
- {\cal L}_0 u &=& \lambda_{\Omega} u & {\rm in} \ \ \Omega,\\
u &=& 0  & {\rm on} \ \ \partial\Omega
\end{array}
\right.
\end{equation}

\begin{theorem} \label{teorema2}
Let $\Omega_0$ be a bounded open subset of $\R^n$ satisfying the divergence theorem for an unit normal vector field $\nu$ defined on $\partial \Omega_0$ in almost everywhere and oriented outwards $\Omega_0$. Let $\varphi: \R \to E^1(\overline{\Omega}_0)$ be a $C^{r+1}$ curve, $r \geq 1$, such that $\varphi(0) = 1_{\Omega_0}$. Denote $\overline \varphi = \varphi'(0)$. Let $\lambda_0$ be a multiple eigenvalue of the operator ${\cal L}_0$ on $\Omega_0$ and of its adjoint ${\cal L}^*_0$ with eigenspaces of same dimension being generated by, respectively, $\{u_0, u_1, \ldots, u_k\}$ and $\{v_0, v_1, \ldots, v_k\}$ in $H^1_0(\Omega_0)$ with $k \geq 1$. Consider the $(k+1) \times (k+1)$ real matrices $A$ and $B$ with elements, respectively,

\[
A_{ij} = (u_i, v_j)_{\overline\varphi}: = \int_{\partial \Omega_0} a(x) \frac{\partial
u_i}{\partial \nu} \frac{\partial v_j}{\partial \nu}\; \o \phi \cdot \nu\; dS
\]

\n and

\[
B_{ij} = \langle u_i, v_j \rangle_{L^2} := \int_{\Omega_0} u_i v_j\; dx\, ,
\]
where $a(x) = \sum_{i,j} a_{ij}(x) \nu_i \nu_j$. Assume that $\lambda_0$ has algebraic and geometric multiplicities equal to $k$ and that the polynomial $\chi(s) = \det(A - sB)$ has a simple real zero $\mu$. Assume also that the eigenfunctions $u_i$ and $v_i$ are of $C^1$ class up to the boundary. Then, there exist $C^{r}$ functions $t \in (-\varepsilon, \varepsilon) \mapsto \lambda(t)$ and $t \in (-\varepsilon, \varepsilon) \mapsto u(t)$ such that $\lambda(0) = \lambda_0$, $u(0) \in \langle u_0, u_1, \ldots, u_k \rangle$ and $(\lambda(t), u(t)) \in \Lambda(-{\cal L}_0, \varphi_t(\Omega_0))$ for $t \in (-\varepsilon, \varepsilon)$. Moreover, we have $\lambda'(0) = -\mu$.
\end{theorem}

\n {\bf Proof of Theorem 3.} Denote $X = \R \times H^1_0(\Omega_0)$ and $Y = \R \times H^{-1}(\Omega_0)$. Let $f: E^1(\overline{\Omega}_0) \times X \to Y$ be the function given by

\[
f(\phi,\lambda, u) = \left(\frac 12 \int_{\Omega_0} u^2 dx , T^\phi_{\Omega_0}(u) - \lambda I^\phi_{\Omega_0}(u) \right)\; ,
\]

\n where $T^\phi_{\Omega_0}(u)(v) = B_{\phi(\Omega_0)}(\phi^* u, \phi^* v)$, with $\phi^* u (x) = u(\phi^{-1}(x))$ and

\[
B_{\phi(\Omega_0)}(u,v) = \int_{\phi(\Omega_0)} \sum_{i,j=1}^n a_{ij}(x) \partial_i u \partial_j v + \sum_{i=1}^n b_i(x) \partial_i u v + c(x) u v\; dx\, ,
\]

\n and $I^\phi_{\Omega_0}(u)(v) = J_{\phi(\Omega_0)}(\phi^* u, \phi^* v)$, with

\[
J_{\phi(\Omega_0)}(u, v) = \int_{\phi(\Omega_0)} u v\; dx\, .
\]

\n Note that $T^\phi_{\Omega_0}(u)$ and $I^\phi_{\Omega_0}(u)$ are linear and continuous, so that $f$ is well-defined.

Let $w \in \R^{k+1} \setminus \{0\}$ be such that $A^Tw = \mu B^Tw$ and consider $u^* = \sum w_i u_i$. Consider the function $F: \R \times X \to Y$ defined by $F(t, \lambda, u) = f(\varphi(t), \lambda, u)$. We now prove that $F$ is in the conditions of Theorem \ref{DIft} at the point $(\lambda_0, u^*)$.

Firstly, $F$ is of $C^{r+1}$ class by Proposition 2.1 of \cite{HM}. Set $K = \Ker({\cal L}_0 + \lambda_0) \subset H^1_0(\Omega_0)$ and $L = \Ker({\cal L}^*_0 + \lambda_0) \subseteq H^1_0(\Omega_0)$. It is clear that
\[
\Img(D_X F(0, \lambda_0, u^*)) = \R \times (\langle u^* \rangle \oplus L^\bot) \subseteq \R \times H^{-1}(\Omega_0)
\]

\n and

\[
\Ker(D_X F(0, \lambda_0, u^*)) = \{0\} \times (K \cap \langle u^* \rangle^\bot) \subseteq \R \times H^1_0(\Omega_0)\, ,
\]

\n where $M^\bot$ denotes the orthogonal in $H^{-1}(\Omega_0)$ of the subspace $M$ of $H^1_0(\Omega_0)$.

Notice that $F(0, \lambda_0, u) = (\frac 12 , 0)$ for all $u \in K$ with $\|u\|_{L^2} = 1$. Then, $F(0,\cdot)$ is constant in a $k$-dimensional sphere with tangent space at $u^*$ equal to $\{0\} \times (K \cap \langle u^*\rangle ^\bot) = \Ker(D_X F(0, \lambda_0, u^*)$. Therefore, the conditions (a) and (b) of Theorem \ref{DIft} are fulfilled.
	
By direct computation, we have

\[
F_t(0, \lambda_0, u) = (0, -({\cal L}_0 + \lambda_0)\frac{\partial u}{\partial \o\varphi})
\]

\n for all $u \in K$. On the other hand, the divergence theorem provides

\[
\langle -({\cal L}_0 + \lambda_0)\frac{\partial u}{\partial \o\varphi}, v \rangle_{L^2} = \int_{\partial \Omega_0} a(x) \frac{\partial u}{\partial \nu} \frac{\partial v}{\partial \nu} \o\varphi \cdot \nu\; ds
\]

\n for all $u \in \Ker({\cal L}_0 + \lambda_0), v \in \Ker({\cal L}_0^* + \lambda_0)$. So, writing $F_t = (F_t^1, F_t^2)$, for any $v \in L$, we derive
	
\[
\langle F^2_t(0, \lambda_0, u), v \rangle_{L^2} = (u, v)_{\o\varphi}\, .
\]
	
Let $\tilde{u} = F^2_t(0, \lambda_0, u^*) = -({\cal L}_0 + \lambda_0)\frac{\partial u^*}{\partial \o\varphi}$. By the definition of $u^*$, for each $j$, we have
	
\[
\langle \tilde{u}, v_j \rangle_{L^2} = (u^*, v_j)_{\o\varphi} = \sum A_{ij} w_i = \mu \sum B_{ij} w_i = \mu \langle u^*, v_j\rangle_{L^2}\, ,
\]

\n meaning that $\tilde{u} - \mu u^* \bot v$ for all $v \in L$ and $\tilde{u} \in \langle u^* \rangle \oplus L^\bot$.
	
Then,

\[
F_t(0, \lambda_0, u^*) = (0, \tilde{u}) \in \Img(D_X F(0, \lambda_0, u^*))\, ,
\]

\n so $F_t(0, \lambda_0, u^*) = D_X F(0,\lambda_0, u^*)[\o \lambda, \o u]$ for some $\o \lambda \in \R$ and $\o u \in H^1_0(\Omega_0)$, and actually $\o \lambda = \mu$. Thus, the condition (c) of Theorem \ref{DIft} is satisfied.

In order to verify the condition (d) of Theorem \ref{DIft}, we first write the directional derivative
\[
F_{(\mu, \o v)}(0, \lambda, u) = D_X F(0, \lambda, u)[(\mu, \o v)] = ( \langle \o v, u \rangle_{L^2}, ({\cal L}_0 + \lambda)(\o v) + \mu u)
\]
	
\n Thus, for any $\o u \in K$,
	
\[
D(F_t - F_{(\mu, \o v)})(0, \lambda_0, u^*)[0 , \o u] = (-\langle \o v, \o u \rangle_{L^2} , -({\cal L}_0 + \lambda_0)\frac{\partial \o u}{\partial \o\varphi} - \mu \o u)
\]
	
Notice now that

\[
\coker(DF(0,\lambda_0, u^*) = \frac {\R \times (H^1_0(\Omega_0))^*}{\R \times (\langle u^* \rangle \oplus L^\bot)} = \{0\} \times (L \cap \langle u^* \rangle^\bot)^*
\]
	
\n where the last identification is the natural one.
	
Then we need to check that the bilinear function
	
\[
G:(K \cap \langle u^* \rangle^\bot) \times (L \cap \langle u^* \rangle^\bot) \to \R
\]
\[
G(\o u, \o v) = \langle D(F^2_t - F^2_{(\mu, \o x)})(0, \lambda_0, u^*)[(0, \o u)], \o v \rangle = (\o u, \o v)_{\o\phi} - \mu \langle\o u, \o v\rangle
\]

\n is non-degenerate. To this end, take a basis of $K$ of the form $\{u_0^*, u_1^*, \ldots, u_k^*\}$, where $u_0^* = u^*$ and $u_i^* \bot u^*$ for $i = 1, \ldots, k$,
and a basis of $L$ of the form $\{v_0^*, v_1^*, \ldots, v_k^*\}$, where $v_i^* \bot u^*$ for $i = 1, \ldots, k$.
Since the algebraic and geometric multiplicities of $\lambda_0$ are equal, we have $u^* \not\in L^\bot$, so $\dim(K \cap \langle u^* \rangle^\bot) = \dim(\langle u^* \rangle \oplus L^\bot) = k$.
For $\nu \in \R$, consider the $(k+1) \times (k+1)$ matrix
\[
M(\nu)_{ij} = (u_i^*, v_j^*)_{\o\phi} - \nu \langle u_i^*, v_j^*\rangle_{L^2} = C^T \cdot (A - \nu B) \cdot D\, ,
\]
where $C$ and $D$ are the matrices of change of basis.
Then
\[
M(\nu) =
\left(
\begin{array}{cc}
(\mu - \nu) \langle u^*, v_0^* \rangle_{L^2} & *\\
0 & N(\nu)
\end{array}
\right)
\]
where $N(\mu)$ is the $k\times k$ matrix of the bilinear form $G$ in the bases $\{u_1^*, \ldots, u_k^*\}$ and $\{v_1^*, \ldots, v_k^*\}$, respectively.
Then, we have
\[
\det(M(\nu)) = \det(C) \det(D) \det(A - \nu B) = \langle u^*, v_0^* \rangle_{L^2}  (\mu - \nu)  \det(N(\nu))
\]
and, by the simplicity assumption of $\mu$, we obtain $\det(N(\mu)) \neq 0$.
	
So, it follows the first conclusion by Theorem \ref{DIft}. Finally, using the fact that $(\lambda'(0), u'(0)) + (\mu, \o u) \in \{0\} \times K$, we deduce the desired formula $\lambda'(0) = -\mu$.\hfill$\Box$

\begin{remark}
The previous proof can be easily adapted to cover the case of complex zeros of $\chi(s)$. Precisely, let $\lambda_0$ be a real eigenvalue of multiplicity $k+1$ of the {\it real} operator $\mathcal L_0$ and assume that the polynomial $\chi(s)$ has a simple complex zero $\mu$. Then, considering the same function $f$, we deduce the existence of a complex branch of eigenvalues $\lambda: (-\epsilon, \epsilon) \to \C$ with $\lambda(0) = \lambda_0$ and $\lambda'(0) = -\mu$.
\end{remark}

\section{Applications}

This last section is devoted to some applications to the $C^r$ dependence problem of Dirichlet eigenvalues with respect to one-parameter $C^{r+1}$ perturbations of the domain with $r \geq 1$. Particularly, we consider the question on disks, squares and balls.

\n {\bf (I) Regularity on disks.} Let $\Omega_0$ be the unit disk centered at the origin in $\R^2$ and $\varphi: \R \to E^1(\overline{\Omega}_0)$ be a $C^{r+1}$ one-parameter perturbation such that $\varphi(0) = 1_{\Omega_0}$. Denote $\overline \varphi = \varphi'(0)$.

The multiple Dirichlet eigenspaces of the Laplace operator on $\Omega_0$ are spanned by functions of the form $u(\rho, \theta) = f(\rho) \sin(k \theta)$ and $v(\rho, \theta) = f(\rho) \cos(k \theta)$ written in polar coordinates. Here $k \in \mathbb N$ and $f(\rho) = J_k(j \rho)$, where $J_k$ is the Bessel of the first kind of order $k$ and $j$ is some zero of $J_k$.

Since $u$ and $v$ are $L^2$ orthogonal, we may assume that $f$ is normalized in such a way that $B$ is the identity matrix.
The normal derivatives are given by
\[
\frac {\partial u}{\partial \nu} = f'(1) \sin(k \theta)
\]
and

\[
\frac {\partial v}{\partial \nu} = f'(1) \cos(k \theta)\, .
\]

Let
\[
\delta (\theta) = \langle \o\phi(\cos \theta,\sin \theta) , (\cos \theta,\sin \theta) \rangle\, .
\]

Let $A$ be the matrix of Theorem \ref{teorema1}. Then the discriminant of its characteristic polynomial is

\[
4 b^2 + a^2  = 4 \left( \int_0^{2\pi} \sin(k \theta) \cos(k \theta) \delta(\theta) d \theta \right)^2 + \left( \int_0^{2\pi} (\cos(k \theta)^2 - \sin(k \theta)^2)\delta(\theta) d \theta \right)^2
\]

\[ = \left( \int_0^{2\pi} \sin(2 k \theta) \delta(\theta) d \theta \right)^2 + \left( \int_0^{2\pi} \cos(2 k \theta) \delta(\theta) d \theta \right)^2  = |\hat\delta(2 k)|^2
\]
which is the $2k$-Fourier coefficient of $\delta$, where $a = |A_{11} - A_{22}|$ and $b = A_{12} = A_{21}$.

Therefore, Theorem \ref{teorema1} provides that all Dirichlet double eigenvalues of the Laplace operator on disks are of $C^{r}$ class close to $t = 0$ whenever $\hat\delta(2k) \neq 0$ for all $k$. On the other hand, the $2k$-Fourier coefficient of $\delta$ to be nonzero for all $k$ is a generic property with respect to perturbations $\varphi$.

In order to illustrate this application, we consider the explicit perturbation $\varphi_t(z) = z + t \sum_{k \geq 0} a_k z^k$, where $z \in \C$, $(ka_k) \in l^1(\R)$ and $a_{2k + 1} \neq 0$ for all $k$. Then, all Dirichlet eigenvalues of the Laplace operator on the unit disk are of $C^\infty$ class close to $t = 0$.\\

\n {\bf (II) Regularity on squares.} Let $\Omega_0 = [0,\pi]^2$ be a square in $\R^2$ and $\varphi: \R \to E^1(\overline{\Omega}_0)$ be a $C^{r+1}$ one-parameter perturbation such that $\varphi(0) = 1_{\Omega_0}$. Denote $\overline \varphi = \varphi'(0)$.

Multiple eigenspaces of $\Omega_0$ with eigenvalue $\lambda$ are spanned by functions of the form
\[
u_\sigma(x_1, x_2) = \frac 2 \pi \sin( \sigma_1 x_1 ) \sin( \sigma_2 x_2 )
\]
with $\sigma_1^2 + \sigma_2^2 = \lambda$ and $\sigma_i \in \N$.

Let us assume for simplicity that the above Diophantine equation for $\sigma_i$ has only two solutions, namely $\sigma = (\sigma_1, \sigma_2)$ and $\o\sigma = (\sigma_2, \sigma_1)$, so the eigenspaces are two-dimensional.

Consider the functions $\eta, \mu:[0, 2\pi] \to \R$

\[
\eta(t) = \left\{\begin{array}{c} -{\o\phi}_2(t,0) \hbox{ if } t<\pi \\ {\o\phi}_2(t, \pi)  \hbox{ if } t \geq \pi \end{array}\right. ,\;\;\;
\mu(t) = \left\{\begin{array}{c} -{\o\phi}_1(0,t) \hbox{ if } t<\pi \\ {\o\phi}_
1(\pi, t)  \hbox{ if } t \geq \pi \end{array}\right.
\]
and the Fourier's coefficients

\[
\hat\eta(k) = \int_0^{2\pi} \eta(t) \cos(k t) dt\;,\;\; \hat\mu(k) = \int_0^{2\pi} \mu(t) \cos(k t) dt\, .
\]
Straightforward computations show that for $\gamma, \delta \in \{\sigma, \o\sigma\}$, we have

\[
A_{\gamma \delta} = \frac{4 \gamma_1 \delta_1}{\pi^2} \int_0^\pi \mu(t) \sin(\gamma_2 t)\sin(\delta_2 t) + (-1)^{\gamma_1+\delta_1} \mu(\pi + t) \sin(\gamma_2 t)\sin(\delta_2 t) dt
\]

\[
\frac{4 \gamma_2 \delta_2}{\pi^2} \int_0^\pi \eta(t) \sin(\gamma_1 t)\sin(\delta_1 t) + (-1)^{\gamma_2+\delta_2} \eta(\pi + t) \sin(\gamma_1 t)\sin(\delta_1 t)dt\, .
\]
and noting that $(-1)^{\gamma_1+\delta_1} = (-1)^{\gamma_2+\delta_2}$, we may change variables in the second term of each integrand to obtain

\[
A_{\gamma \delta} = \frac{4 \gamma_1 \delta_1}{\pi^2} \int_0^{2\pi} \mu(t) \sin(\gamma_2 t)\sin(\delta_2 t) dt +\frac{4 \gamma_2 \delta_2}{\pi^2} \int_0^{2\pi} \eta(t) \sin(\gamma_1 t)\sin(\delta_1 t) dt\, .
\]

Now using standard trigonometric identities, we have

\[
A_{\gamma \delta} = \frac{2 \gamma_1 \delta_1}{\pi^2} (\hat\mu(\gamma_2-\delta_2) - \hat\mu(\gamma_2+\delta_2)) + \frac{2 \gamma_2 \delta_2}{\pi^2} (\hat\eta(\gamma_1-\delta_1) - \hat\eta(\gamma_1 + \delta_1))\, .
\]
Then the coefficients of $A$ are

\[
A_{\sigma \sigma} = \frac {2 \sigma_1^2}{\pi^2}(\hat\mu(0) - \hat\mu(2 \sigma_2))+  \frac {2 \sigma_2^2}{\pi^2}(\hat\eta(0) - \hat\eta(2 \sigma_1)))\, ,
\]

\[
A_{\o\sigma \o\sigma} = \frac {2 \sigma_2^2}{\pi^2}(\hat\mu(0) - \hat\mu(2 \sigma_1))+  \frac {2 \sigma_1^2}{\pi^2}(\hat\eta(0) - \hat\eta(2 \sigma_2)))\, ,
\]

\[
A_{\sigma \o\sigma} = A_{\sigma \o\sigma} = \frac {2 \sigma_1 \sigma_2}{\pi^2}(\hat\mu(\sigma_1-\sigma_2) - \hat\mu(\sigma_1+\sigma_2) + \hat\eta(\sigma_1-\sigma_2) - \hat\eta(\sigma_1+\sigma_2))\, .
\]

Finally, $A$ has simple eigenvalues if and only if either $A_{\sigma \sigma} \neq A_{\o\sigma \o\sigma}$ or $A_{\sigma \o\sigma} \neq 0$. Or equivalently, if and only if one of the following conditions occur:

\[
(\hat\eta+\hat\mu) (\sigma_1+\sigma_2) \neq (\hat\eta+\hat\mu)(\sigma_1-\sigma_2)\, ,
\]

\[
\sigma_1^2(\hat\mu-\hat\eta)(2\sigma_2) - \sigma_2^2(\hat\mu-\hat\eta)(2\sigma_1) \neq (\sigma_1^2-\sigma_2^2) (\hat\mu-\hat\eta)(0)\, .
\]
Thus, since these two situations are generic with respect to perturbations $\varphi$, Theorem \ref{teorema1} provides that all Dirichlet eigenvalues of the Laplace operator on squares are of $C^{r}$ class close to $t = 0$ for almost every perturbation $\varphi$.
\\

\n {\bf (III) Regularity on balls in $\R^n$.} Let $\Omega_0$ be the unit ball centered at the origin in $\R^n$ and $\varphi: \R \to E^1(\overline{\Omega}_0)$ be a $C^{r+1}$ one-parameter perturbation such that $\varphi(0) = 1_{\Omega_0}$. Denote $\overline \varphi = \varphi'(0)$.

The first multiple eigenspace in the unit ball of $\R^n$ is spanned by functions of the form $u_i(x) = f(\rho) \frac {x_i} \rho$, where $\rho = \|x\|$, $i = 1, \ldots, n$.
The derivatives in the normal direction coincide in $\s^{n-1}$ with $x_1, \ldots, x_n$, except for a constant.

Denote
\[
\delta(\theta) = \langle \o\phi(\theta),\nu(\theta)\rangle\, ,
\]
where $\nu$ is the unit normal vector field of $\s^{n-1}$ oriented outwards. Its Fourier coefficients of degree $2$ are

\[
F_i = \int_{\s^{n-1}} (\theta_i^2 - \theta_1^2) \delta(\theta) dS,\, i = 2, \ldots, n\, ,
\]

\[
C_{ij} = \int_{\s^{n-1}} \theta_i \theta_j \delta(\theta) dS,\; i,j = 1, \ldots, n,\, i > j\, .
\]

Let $a = \int_{\s^{n-1}} \theta_1^2 \delta(\theta) dS$, we compute

\[
A - a Id =
\left(
\begin{array}{ccccccc}
0 & C_{21} & C_{31} & \ldots & C_{n1}\\
C_{21} & F_2 & C_{32} &\ldots & C_{n2}\\
C_{31} & C_{32} &F_3 &\ldots & C_{n2}\\
\vdots\\
C_{n1} & C_{n2} & C_{n3} &\ldots & F_n\\

\end{array}
\right).
\]
Since $\chi_A(s) = \chi_{A-a Id}(s-a)$, the simplicity of the eigenvalues of $A$ can be translated from the second order Fourier coefficients of $\delta$. In particular, the characteristic polynomial $\chi_A$ of the matrix $A$ is constructed from sums and products of these latter ones. Thus, the discriminant of $\chi_A$ is nonzero for almost every perturbation $\varphi$, so that all $n$ eigenvalues are of $C^{r}$ class close to $t = 0$.

The problem of determining the simplicity of the eigenvalues of $A$ in general seems to be a complicated one for balls in $\R^n$ with $n \geq 3$ and we don't hope to find a simple solution in these cases.\\

\n {\bf (IV) Two disjoint domains.} Let us apply Theorem \ref{teorema1} into a quite situation.

Consider $C^{r+1}$ one-parameter perturbations $\varphi: \R \to E^1(\overline{U}_0)$ and $\psi: \R \to E^1(\overline{V}_0)$. Denote $U_t = \varphi_t(U_0)$ and $V_t = \psi_t(V_0)$ with principal eigenvalues $\mu_t$ and $\nu_t$, respectively, where $U_0$ and $V_0$ are disjoint bounded open subsets of $\R^n$ but equal after a translation.

Let $W = U_0 \cup V_0$. Since $\mu_0 = \nu_0$, the eigenspace is spanned by functions $u_t$ and $v_t$ with support at $U_t$ and $V_t$, respectively.
Assuming that $\|u_0\|_{L^2} = \|v_0\|_{L^2} = 1$, the inner product matrix $B$ is the identity and $A$ is diagonal.
Moreover, $A_{11} = -\mu'(0)$ and $A_{22} = -\nu'(0)$ by the Hadamard formula, see for example \cite{HM}.

In this case we evoke Theorem \ref{teorema1} which states that $\mu_t$ and $\nu_t$ are $C^{r}$ functions on $t$, provided that $A_{11} \neq A_{22}$.
However, observe that, unlike $u_t$ and $v_t$, not every eigenfunction of $W_t = U_t \cup V_t$ is smooth on $t$.
For example, if $w_0 = u_0 + v_0$ the local uniqueness of the principal eigenvalues of $U_0$ and $V_0$ guarantees that the unique possible continuation is $w_t = u_t + v_t$ which is not an eigenfunction unless $\mu_t = \nu_t$ for every $t$ small enough. Part of the proof of Theorem \ref{teorema1} is to find the eigenfunction having a continuation.\\

\n {\bf Acknowledgments:} The first author was supported by CAPES (BJT 064/2013) and Fapemig. The second author acknowledges the support provided by CAPES (BEX 6961/14-2), CNPq (PQ 306406/2013-6) and Fapemig (PPM 00223-13).

\end{document}